\documentclass[12pt]{article}
\usepackage{amsmath,amssymb,amsthm,bm,comment}
\usepackage{tikz}   
\usetikzlibrary{arrows.meta,calc}
\usepackage{array}
\usepackage{graphicx}

\def\beq{\begin{equation}}
\def\eqn#1{\beq\label{#1}}
\def\ee{\end{equation}}

\def\bb {\begin {eqnarray}}
\def\eqnn#1{\bb\label{#1}}
\def\eea {\end {eqnarray}}

\newcommand{\eqna}[1]{\begin{subequations} \label{#1}
\begin{eqnarray}}
\def\eena{\end{eqnarray}
\end{subequations}}

\def\Lp{\L^+}

\def\chh{\chi^+}

\def\white#1{\mathop{\bigcirc}\limits_{#1}}

\def\blackr#1{\overset{#1}{\mathop{\Bbullet}}}

\def\nn{\nonumber}

\def\llr{\longrightarrow}
\def\({\left(}
\def\){\right)}
\def\eps{\epsilon}

\def\ha{{\textstyle{1\over2}}}

  \def\tV{{\tilde V}}

\def\r{\rho}
\def\cgc{{\cg^\bbc}}

\def\a{\alpha}
\def\b{\beta}

\def\vr{\vert}
\def\g{\gamma}

\def\y{\eta}

\def\D{{\Delta}}

\def\bbr{{I\!\!R}}
\def\bbn{I\!\!N}
\def\bbc{{C\kern-8.5pt I}\,}
\def\bac{\bbc} 
\def\bab{{C\kern-4.5pt I}}

\def\bbc{\mathbb{C}}

\def\bbr{\mathbb{R}}
\def\bbn{\mathbb{N}}

\def\L{\Lambda}

\def\rank{{\rm rank}}

\def\riga{-\kern-4pt - \kern-4pt -}
\font\fat=cmsy10 scaled\magstep5

\def\Bbullet{\raise-3pt\hbox{\fat\char"0F}}

\def\ca{{\cal A}}  \def\cc{{\cal C}}
\def\cd{{\cal D}}  \def\cf{{\cal F}}
\def\cg{{\cal G}} \def\ch{{\cal H}} \def\ci{{\cal I}}
 \def\ck{{\cal K}} 
\def\cm{{\cal M}} \def\cn{{\cal N}} 
\def\cp{{\cal P}}  
 \def\ct{{\cal T}}

\def\ido{intertwining differential operator}
\def\idos{intertwining differential operators}

 \def\ha{{\textstyle{\frac{1}{2}}}}

\def\eps{\epsilon}

\def\ca{{\cal A}}
\def\nn{\nonumber}

\input epsf.tex
\newcount\figno
\def\fig#1#2#3{
\par\begingroup\parindent=0pt\leftskip=1cm\rightskip=1cm\parindent=0pt
\baselineskip=11pt \global\advance\figno by 1 
\epsfxsize=#3 \centerline{\epsfbox{#2}} \vskip 12pt
#1\par
\endgroup\par}
\def\figlabel#1{\xdef#1{\the\figno}}
\def\encadremath#1{\vbox{\hrule\hbox{\vrule\kern8pt\vbox{\kern8pt
\hbox{$\displaystyle #1$}\kern8pt} \kern8pt\vrule}\hrule}}

\begin{document}

\begin{center}

{\Large{Invariant Differential Operators for Non-Compact}}
\vskip 0.2cm
{\Large{ Lie Groups: the   $Sp(n,1)$  Case}}

\vspace{10mm}

{{\large\bf  N. Aizawa$^1$, ~V.K. Dobrev$^2$}
	\\[10pt]
$^1$Department of Physics, Graduate School of Science,\\ Osaka Metropolitan University,   \\
Nakamozu Campus, Sakai, Osaka 599-8531, Japan
\\[10pt]
$^2$Institute of Nuclear Research and Nuclear Energy,\\
 Bulgarian Academy of Sciences, \\
72 Tsarigradsko Chaussee, 1784 Sofia, Bulgaria}

\end{center}

\vspace{10mm}

\vskip 0.8cm

\begin{abstract}In the present paper we continue the project of systematic
construction of invariant differential operators on the example of
the non-compact  algebras  $sp(n,1)$. Our
choice of these algebras is motivated by the fact that   they belong
to a narrow class of algebras, which are of split rank one, of which class
the other cases were studied, some long time ago. We concentrate on the case $n=2$.   
We give the main multiplets and the main  reduced multiplets of indecomposable elementary representations for, including the necessary data for all relevant invariant differential
operators. We also present explicit expressions for the singular vectors and the \idos.

\end{abstract}

\vfill\eject

\section{Introduction}

Consider a Lie group ~$G$, e.g., the Lorentz, Poincar\'e, conformal
groups, and differential equations
$$ \ci~f ~~=~~ j $$
which are $G$-invariant. These
play a very important role in the description of physical
symmetries - recall, e.g., the early examples of Dirac, Maxwell, d'Allembert,
equations and nowadays the latest applications
of (super-)differen\-tial operators in conformal field theory,
supergravity, string theory. Naturally, it is important to construct systematically
such invariant equations and operators.

In the present paper we  focus on the groups ~$Sp(n,1)$, ($n\geq 2$).
which are very interesting   since they belong
to a narrow class of algebras, which are of split rank one, of which class
the other cases were studied, some long time ago. Namely, the Euclidean conformal group $SO(N,1)$
was studied in \cite{DMPPT}. By a different approach in \cite{Zhel} were studied the cases
$SO(N,1)$, $SU(N,1)$. Recently, in \cite{Dobf4} was studied the case ~$F''_4$.

$Sp(n,1)$ also has some physical relevance. 
The importance of symplectic groups $Sp(n,\mathbb{R})$ in physics is very well-known. 
In any physical applications of groups, the knowledge of their representations is fundamental and necessary. 
One may use the representations of $Sp(n,1)$ in order to investigate those of $Sp(n+1,\mathbb{R}).$ 
Taking a quotient of $Sp(n,1)$ by an appropriate subgroup, one may have the quaternionic projective space on which Laplace-Beltrami operator, which is also physically important, is defined, e.g. \cite{Pajas}. Such projective space is also used to define a model of supergravity with higher derivative couplings \cite{SUGRA}.

The present paper is organized as follows. In section 2 we give the
preliminaries recalling and adapting facts from
\cite{Dob}. In Section 3 we specialize to the ~$sp(n,1)$~
case, concentrating on the case $n=2$. In Section 4 we present the results on the multiplet
classification of the representations, the singular vectors and the corresponding \idos.

\section{Preliminaries} \label{SEC:Pre}

The text in this section is standard, taken, e.g., from  \cite{Dob}. Let $G$ be a semisimple non-compact Lie group, and $K$ a
maximal compact subgroup of $G$. Then we have an Iwasawa
decomposition ~$G=KA_0N_0$, where ~$A_0$~ is abelian simply
connected vector subgroup of ~$G$, ~$N_0$~ is a nilpotent simply
connected subgroup of ~$G$~ preserved by the action of ~$A_0$.
Further, let $M_0$ be the centralizer of $A_0$ in $K$. Then the
subgroup ~$P_0 ~=~ M_0 A_0 N_0$~ is a minimal parabolic subgroup of
$G$. A parabolic subgroup ~$P ~=~ M' A' N'$~ is any subgroup of $G$
(including $G$ itself) which contains a minimal parabolic subgroup.

The importance of the parabolic subgroups comes from the fact that
the representations induced from them generate all (admissible)
irreducible representations of $G$ \cite{Lan}. For the
classification of all irreducible representations it is enough to
use only the so-called {\it cuspidal} parabolic subgroups
~$P=M'A'N'$, singled out by the condition that ~rank$\, M' =$
rank$\, M'\cap K$ \cite{Zhea,KnZu}, so that $M'$ has discrete series
representations \cite{Har}. However, often induction from
non-cuspidal parabolics is also convenient, cf.
\cite{EHW,Dob}.

Let ~$\nu$~ be a (non-unitary) character of ~$A'$, ~$\nu\in\ca'^*$,
let ~$\mu$~ fix an irreducible representation ~$D^\mu$~ of ~$M'$~ on
a vector space ~$V_\mu\,$.

 We call the induced
representation ~$\chi =$ Ind$^G_{P}(\mu\otimes\nu \otimes 1)$~ an
~{\it elementary representation} (ER) of $G$ \cite{DMPPT}. (These are
called {\it generalized principal series representations} (or {\it
limits thereof}) in \cite{Knapp}.) Their spaces of functions are:
\eqn{fun} \cc_\chi ~=~ \{ \cf \in C^\infty(G,V_\mu) ~ \vr ~ \cf
(gman) ~=~ e^{-\nu(H)} \cdot D^\mu(m^{-1})\, \cf (g) \} \ee where
~$a= \exp(H)\in A'$, ~$H\in\ca'\,$, ~$m\in M'$, ~$n\in N'$. The
representation action is the $left$ regular action: \eqn{lrr}
(\ct^\chi(g)\cf) (g') ~=~ \cf (g^{-1}g') ~, \quad g,g'\in G\ .\ee

For our purposes we need to restrict to ~{\it maximal}~ parabolic
subgroups ~$P$, (so that $\rank\,A'=1$), that may not be cuspidal.
For the representations that we consider the character ~$\nu$~ is
parameterized by a real number ~$d$, called the conformal weight or
energy.

Further, let ~$\mu$~ fix a discrete series representation ~$D^\mu$~
of $M'$ on the Hilbert space ~$V_\mu\,$, or the so-called limit of a
discrete series representation (cf. \cite{Knapp}). Actually, instead
of the discrete series we can use the finite-dimensional
(non-unitary) representation of $M'$ with the same Casimirs.

An important ingredient in our considerations are the ~{\it
highest/lowest weight representations}~ of ~$\cg$. These can be
realized as (factor-modules of) Verma modules ~$V^\L$~ over
~$\cg^\bac$, where ~$\L\in (\ch^\bac)^*$, ~$\ch^\bac$ is a Cartan
subalgebra of ~$\cg^\bac$, weight ~$\L = \L(\chi)$~ is determined
uniquely from $\chi$ \cite{Dob}. In this setting we can consider
also unitarity, which here means positivity w.r.t. the Shapovalov
form in which the conjugation is the one singling out $\cg$ from
$\cg^\bac$.

Actually, since our ERs may be induced from finite-dimensional
representations of ~$\cm'$~ (or their limits) the Verma modules are
always reducible. Thus, it is more convenient to use ~{\it
generalized Verma modules} ~$\tV^\L$~ such that the role of the
highest/lowest weight vector $v_0$ is taken by the
(finite-dimensional) space ~$V_\mu\,v_0\,$. For the generalized
Verma modules (GVMs) the reducibility is controlled only by the
value of the conformal weight $d$. Relatedly, for the \idos{} only
the reducibility w.r.t. non-compact roots is essential.

One main ingredient of our approach is as follows. We group the
(reducible) ERs with the same Casimirs in sets called ~{\it
multiplets} \cite{Dob,Dobmul}. The multiplet corresponding to fixed
values of the Casimirs may be depicted as a connected graph, the
vertices of which correspond to the reducible ERs and the lines
between the vertices correspond to intertwining operators. The
explicit parametrization of the multiplets and of their ERs is
important for understanding of the situation.

In fact, the multiplets contain explicitly all the data necessary to
construct the \idos{}. Actually, the data for each \ido{} consists
of the pair ~$(\b,m)$, where $\b$ is a (non-compact) positive root
of ~$\cg^\bac$, ~$m\in\bbn$, such that the BGG \cite{BGG} Verma
module reducibility condition (for highest weight modules) is
fulfilled: \eqn{bggr} (\L+\r, \b^\vee ) ~=~ m \ , \quad \b^\vee
\equiv 2 \b /(\b,\b) \ \ee 
where $\rho$ is the half-sum of positive roots. 
When (\ref{bggr}) holds then the Verma
module with shifted weight ~$V^{\L-m\b}$ (or ~$\tV^{\L-m\b}$ ~ for
GVM and $\b$ non-compact) is embedded in the Verma module ~$V^{\L}$
(or ~$\tV^{\L}$). This embedding is realized by a singular vector
~$v_s$~ determined by a polynomial ~$\cp_{m,\b}(\cg^{\bac}_-)$~ in the
universal enveloping algebra ~$U(\cg^{\bac}_-)\,$, ~$\cg^{\bac}_-$~ is the
subalgebra of ~$\cg^\bac$ generated by the negative root generators
\cite{Dix}.
 More explicitly, \cite{Dob}, ~$v^s_{m,\b} = \cp_{m,\b}\, v_0$ (or ~$v^s_{m,\b} =
 \cp_{m,\b}\, V_\mu\,v_0$ for GVMs).
Then there exists \cite{Dob} an \ido{} \eqn{lido}
\cd_{m,\b} ~:~ \cc_{\chi(\L)} ~\llr ~ \cc_{\chi(\L-m\b)} \ee given
explicitly by: \eqn{invop}\cd_{m,\b} ~=~ \cp_{m,\b}(\widehat{\cg^{\bac}_-})
\ee where ~$\widehat{\cg^{\bac}_-}$~ denotes the $right$ action on the
functions ~$\cf$ which is defined by
	\begin{equation}
	\widehat{X} \mathcal{F}(g) =  \left.\frac{d}{dt} \mathcal{F}(g e^{tX}) \right|_{t=0}, \quad X \in \mathcal{G}^{\mathbb{C}}
\end{equation}		
One may show using the right covariance of $\mathcal{F}$, 
, cf. (\ref{fun}), that the $\mathcal{F}$ has the property of the highest weight \cite{Dob}
\bb
 \widehat{X} \mathcal{F} &=& \Lambda \mathcal{F}, \quad X \in \mathcal{H}^{\mathbb{C}}
 \nonumber \\
 \widehat{X} \mathcal{F} &=& 0, \quad X \in \mathcal{G}^{\mathbb{C}}_+
\eea
where $\cg^{\bac}_+$ is the subalgebra of ~$\cg^\bac$ generated by the negative root generators.

\section{The Non-Compact Lie Algebras $sp(n,1)$}

We start with the complexification $\cg^{\mathbb{C}} = sp(n+1, \mathbb{C})$:
\begin{eqnarray}
	sp(N,\mathbb{C}) = \{   X \in gl(2N,\mathbb{C}) \ | \ {}^tX J_N + J_N X = 0 \},
	\quad
	J_N =
	\begin{pmatrix}
		0 & 1_N \\ -1_N & 0
	\end{pmatrix}
\end{eqnarray}
Writing $ X = \begin{pmatrix}
	A & B \\ C & D
   \end{pmatrix}$
where $ A, B, C, D $ are $N\times N$ complex matrices.
It follows from
\begin{eqnarray}
	{}^tX J_N + J_N X =
	\begin{pmatrix}
		-{}^tC + C & {}^t A + D \\
		-{}^t D - A & {}^t B - B
	\end{pmatrix}
    = 0
\end{eqnarray}
that
\begin{eqnarray}
	X \in sp(N,\mathbb{C}) \ \Rightarrow \
	X =
	\begin{pmatrix}
		A & B \\ C & -{}^t A
	\end{pmatrix},
    \quad
     {}^t B = B, \ {}^t C = C
\end{eqnarray}

Next we need the rank 1 real form of $\cg^\bac$: 
\begin{align}
	  sp(n,1) &= \{  X \in sp(n+1,\mathbb{C}) \ | \ X^{\dagger} y_2 + y_2 X = 0 \},
	  \label{spn1defrel} \\
	  y_2 &
	  =
	  \begin{pmatrix}
	  	\beta_2 & 0 \\ 0 & \beta_2
	  \end{pmatrix},
      \qquad
      \beta_2 =
      \begin{pmatrix}
      	1_{n-1} & 0 & 0
      	\\
      	0 & 0 & 1
      	\\
      	0 & 1 & 0
      \end{pmatrix}.
\end{align}
The matrix $\beta_2$ is related to $\beta_0$ by the unitary transformation
\begin{equation}
	\beta_0 =
	\begin{pmatrix}
		1_{n} & 0 \\ 0 & -1
	\end{pmatrix}
    \mapsto \beta_2 = U \beta_0 U^{\dagger},
    \quad
    U = \frac{1}{\sqrt{2}}
      \begin{pmatrix}
      	\sqrt{2}\cdot 1_{n-1} & 0 & 0
      	\\
      	0 & 1 & -1
      	\\
      	0 & 1 & 1
      \end{pmatrix}.
\end{equation}

The definition \eqref{spn1defrel} constraints  the entries of $A=(a_{ij}), B=(b_{ij}), C=(c_{ij})$ as follows:
\begin{alignat}{3}
  \bar{a}_{ij} + a_{ji} &= 0, &
  \bar{c}_{ij} + b_{ji} &= 0, & \qquad
  & 1 \leq i, j \leq n-1
  \nonumber \\
  \bar{a}_{in} + a_{n+1,i} &= 0, &
  \bar{c}_{in} + b_{n+1,i} &= 0, &
  & 1 \leq i \leq n+1
  \nonumber \\
  \bar{a}_{i,n+1} + a_{n+1,i} &= 0, & \qquad
  \bar{c}_{i,n+1} + b_{n+1,i} &= 0, &
  & 1 \leq i \leq n+1 \label{spn1-1}
\end{alignat}
From the relations one may observe that the following elements are pure imaginary:
\begin{eqnarray}
	 a_{ii}, \ (1 \leq i \leq n-1), \qquad a_{n,n+1}, \qquad a_{n+1,n}   \label{spn1-2}
\end{eqnarray}

 Let  further ~$n\geq 2$.\footnote{The case $n=1$ is not typical for $sp(n,1)$ since $sp(1,1)\cong so(4,1)$
 and the latter group is the Euclidean conformal group of $\bbr^3$ (called also de Sitter group) and the family of $SO(n,1)$ was studied in detail long time ago (cf. \cite{DMPPT}, especially the case $n=4$ in Chapter 7.C).}
 Let ~$\cg ~=~ sp(n,1)$, the  real rank one form
of ~$sp(n+1,\bbc)=\cg^\bac$. The maximal compact subgroup of ~$\cg$~
is ~$\ck \cong sp(n) \oplus sp(1)$, ~$\dim_\bbr\,\cp = 4n$,
~$\dim_\bbr\,\cn = 4n-1$. This algebra has discrete series
representations.

The split rank is equal to $1$, while ~$\cm = sp(n-1) \oplus sp(1)$.

The Satake diagram \cite{Sata} of ~$sp(n,1)$~ is:
\eqn{satdyn}
\blackr{{1}}
\riga
\white{{2}} \riga
  \blackr{{3}}~ \riga \cdots \riga \blackr{{n}}  \Longleftarrow   \blackr{{n+1}}
\ee
where by standard convention the first black dot  represents the $sp(1)$ subalgebra of ~$\cm$, while the
other black dots represent the subalgebra $sp(n-1)$.

We choose a ~{\it minimal/maximal} parabolic ~$\cp=\cm\ca\cn$.~

We label   the signature of the ERs of $\cg$   as follows:
\eqn{sgnd}  \chi ~=~ \{\, m_1, \ldots, m_{n+1} \,\} =
[\, n\, ;\, n_1\,, \ldots,\, n_{n-1}\, ;\, c\, ]  \ ,
\ n,n_j \in \bbn\ , \ c = d-n+\frac{1}{2} \ee
where ~$m_j$~ are the Dynkin labels of $\cgc$,  the entry ~$c$~ fixes the character of $\ca\,$, and the
 entries $n$, $n_j$ are labels of the finite-dimensional  irreps
of $sp(1)\,$, $sp(n-1)$, resp.

Further, we need more explicitly the root system of the algebra
~$sp(n+1,\bbc)$.

In terms of the orthonormal basis $\eps_i\,$, ~$i=1,\ldots,n+1$, the
 positive roots are given by \eqn{spnrpos} \D^+ = \{ \eps_i \pm \eps_j, ~1 \leq
i <j \leq n+1; ~2\eps_i, 1 \leq i \leq n+1\} ,  \ee while the
simple roots are: \eqn{spnrsmp} \pi = \{\g_i = \eps_i - \eps_{i+1},
~1 \leq i \leq n,  ~ \g_{n+1} = 2\eps_{n+1}\} \ee With our choice of
normalization of
  the long roots  ~$2\eps_k$~ have
length 4, while the short roots ~$\eps_i \pm \eps_j$~ have length 2.

From these the $\ck$-compact roots are those that form (by restriction)
the root system of the semisimple part of ~$\ck^\bac$, the rest are
noncompact, i.e., \eqnn{spnrcnc}   \ck{\rm-compact:}&~~~ \a_{ij}
~\equiv~\eps_i  -\eps_j,\ , \quad 1 \leq i < j \leq n+1 \ ,
  \cr  \ck{\rm-noncompact:}&~~~ \b_{ij} ~\equiv~ \eps_i
+\eps_j,\ ,  ~~ 1 \leq i \leq j \leq n+1 \  \eea
Thus, the only
$\ck$-non-compact simple root is ~$\g_n=\b_{n+1,n+1}\,$.

In the above notation the long roots are ~$\b_{jj}$, the rest are short.

We adopt the following ordering of the roots: \eqn{ord} \begin{matrix}
\b_{11}  &&&&& &&&\cr \vee &&&&&&&&\cr \b_{12} &
> & \b_{22} &&&&&\cr \vee && \vee &&&&&&\cr \ldots & \ldots & \ldots
&\ldots & \ldots &&&&\cr \vee  &  & \vee  &  & \ldots & & \vee &&\cr
\b_{1n} &
> & \b_{2n} & > & \ldots & > & \b_{n-1,n} & > &\b_{nn}=\g_n \cr \vee
&  & \vee  &  & \ldots & & \vee &&\cr \a_{1n} & > & \a_{2n} & > &
\ldots & > & \a_{n-1,n} =\g_{n-1}&&\cr \vee  &  & \vee  &  & \ldots
& &  &&\cr \ldots & \ldots & \ldots &\ldots & \ldots &&&&\cr \vee &&
\vee &&&&&&\cr \a_{13} & > & \a_{23}=\g_2 &&&&&\cr \vee && &&&&&&\cr
\a_{12}=\g_1 &  &  &&&&&\cr \end{matrix}\ee This ordering is lexicographical
adopting the ordering of the ~$\eps$~: \eqn{orde} \eps_1 > \cdots >
\eps_n \ee

Further, we shall use the so-called Dynkin labels: \eqn{dynk} m_i
~\equiv~ (\L+\r,\g^\vee_i)  \ , \quad i=1,\ldots,n,\ee where ~$\L =
\L(\chi)$, ~$\r$ is half the sum of the positive roots of
~$\cg^\bac$.

We shall use also   the so-called Harish-Chandra parameters:
\eqn{dynhc} m_\b \equiv (\L+\r, \b^\vee )\ ,\ee where $\b$ is any
positive root of $\cg^\bac$.  In particular, in the
case of the $\ck$-noncompact roots we have: \eqnn{hclab} m_{\b_{ij}} ~&=&~
\Big( \sum_{s=i}^n + \sum_{s=j}^n \Big)   m_s \ , \qquad i<j \ ,\cr
m_{\b_{ii}} ~&=&~  \sum_{s=i}^n    m_s  \eea

There are several types of multiplets: the main type, (which
contains maximal number of ERs/GVMs, the finite-dimensional and the
discrete series representations), and some reduced types of
multiplets.

According to the results of \cite{VKD1} the number of GVMs in a main multiplet should be:
\eqn{numweyl} \frac{ | W(\cg^\bbc) |} {| W(\cm^\bbc)|} ~=~
\frac{ | W(sp(n,1)^\bbc) |} { | W(sp(n-1)^\bbc)| \times |W(sp(1)^\bbc)|}  ~=~ 2n(n+1) \ee

\section{The case ~$sp(2,1)$}

This is the {\bf main topic} of this paper.

Taking into account the relations \eqref{spn1-1} and \eqref{spn1-2}, and taking the minimal parabolic subalgebra, we take the following basis of $sp(2,1)$:
\begin{align}
	{\cal A}_0	  = \mathrm{diag}(0,1,-1,0,-1,1).
\end{align}
The basis of $ \mathcal{M}_0 = sp(1) \oplus sp(1)$
\begin{alignat}{2}
	M_1 &= \mathrm{diag}(i,0,0,-i,0,0), & \quad
	M_2 &= \mathrm{diag}(0,i,i,0,-i,-i),
	\nonumber \\
	M_3  &= \left(
\begin{array}{rrr|rrr}
	& & & 1 & 0 & 0
	\\
	& & & 0 & 0 & 0
	\\
	& & & 0 & 0 & 0
	\\ \hline
	-1 &  0 & 0 & 0 & 0 & 0
	\\
	0 & 0 & 0 & & &
	\\
	0 & 0 & 0 & & &
\end{array}
\right),
&\qquad
M_4 &=
\left(
\begin{array}{rrr|rrr}
	& & & i & 0 & 0
	\\
	& & & 0 & 0 & 0
	\\
	& & & 0 & 0 & 0
	\\ \hline
	i & 0 & 0 & & &
	\\
	0 & 0 & 0 & & &
	\\
	0 & 0 & 0 & & &
\end{array}
\right),
\nonumber \\
     M_5 &=
\left(
\begin{array}{rrr|rrr}
	& & & 0 & 0 & 0
	\\
	& & & 0 & 0 & 1
	\\
	& & & 0 & 1 & 0
	\\ \hline
	0 & 0 & 0 & & &
	\\
	0 & 0 & -1 & & &
	\\
	0 & -1 & 0 & & &
\end{array}
\right),
&
M_6 &=
\left(
\begin{array}{rrr|rrr}
	& & & 0 & 0 & 0
	\\
	& & & 0 & 0 & i
	\\
	& & & 0 & i & 0
	\\ \hline
	0 & 0 & 0 & & &
	\\
	0 & 0 & i & & &
	\\
	0 & i & 0 & & &
\end{array}
\right).
\end{alignat}
The matrices $M_1, M_3, M_4$ span $sp(1)$ and the rest spans the other $sp(1).$
The restricted positive roots of $ \mathcal{A}_0$ are taken to be
$$\begin{matrix}
	L^+_1 &=
	\left(
	\begin{array}{rrr|rrr}
		0 & 0 & 1 & & &
		\\
		-1 & 0 & 0 & & &
		\\
		0 & 0 & 0 & & &
		\\ \hline
		& & & 0 & 1 & 0
		\\
		& & & 0 & 0 & 0
		\\
		& & & -1 & 0 & 0
	\end{array}
	\right),
	& \quad
	L^+_2 &=
	\left(
	\begin{array}{rrr|rrr}
		0 & 0 & i & & &
		\\
		i & 0 & 0 & & &
		\\
		0 & 0 & 0 & & &
		\\ \hline
		& & & 0 & -i & 0
		\\
		& & & 0 & 0 & 0
		\\
		& & & -i & 0 & 0
	\end{array}
	\right),
	\nonumber \\
	L^+_3 &=
	\left(
	\begin{array}{rrr|rrr}
		& & & 0 & 1 & 0
		\\
		& & & 1 & 0 & 0
		\\
		& & & 0 & 0 & 0
		\\ \hline
		0 & 0 & -1 & & &
		\\
		0 & 0 & 0 & & &
		\\
		-1 & 0 & 0 & & &
	\end{array}
	\right),
	&
	L^+_4 &=
	\left(
	\begin{array}{rrr|rrr}
		& & & 0 & i & 0
		\\
		& & & i & 0 & 0
		\\
		& & & 0 & 0 & 0
		\\ \hline
		0 & 0 & i & & &
		\\
		0 & 0 & 0 & & &
		\\
		i & 0 & 0 & & &
	\end{array}
	\right),
	\nonumber \\
	L^+_5 &=
	\left(
	\begin{array}{rrr|rrr}
		0 & 0 & 0 & & &
		\\
		0 & 0 & i & & &
		\\
		0 & 0 & 0 & & &
		\\ \hline
		& & & 0 & 0 & 0
		\\
		& & & 0 & 0 & 0
		\\
		& & & 0 & -i & 0
	\end{array}
	\right),
	&
	L^+_6 &=
	\left(
	\begin{array}{rrr|rrr}
		& & & 0 & 0 & 0
		\\
		& & & 0 & 1 & 0
		\\
		& & & 0 & 0 & 0
		\\ \hline
		0 & 0 & 0 & & &
		\\
		0 & 0 & 0 & & &
		\\
		0 & 0 & -1 & & &
	\end{array}
	\right),
	\nonumber \\
	L^+_7 &=
	\left(
	\begin{array}{rrr|rrr}
		& & & 0 & 0 & 0
		\\
		& & & 0 & i & 0
		\\
		& & & 0 & 0 & 0
		\\ \hline
		0 & 0 & 0 & & &
		\\
		0 & 0 & 0 & & &
		\\
		0 & 0 & i & & &
	\end{array}
	\right)
\end{matrix}$$
and the negative roots are obtained by the Cartan involution:
\begin{equation}
	L_k^- = \theta(L_k^+) :=  y_2 L_k^+ y_2.
\end{equation}
The Cartan subalgebra is spanned by $ \mathcal{A}_0, M_1, M_2. $

\section{Complex Lie alegebra {$sp(3,\mathbb{C})$ }  }

We follow \S 2.5.3 of \cite{VKD1}. Denoting the basis of the Cartan subalgebra of $ sp(3,\mathbb{C})$  by $ H_1, H_2, H_3$,
the dual algebra is given by
\begin{equation}
	\alpha_i(H)_j = \alpha_j(H)_i =
	\begin{cases}
		-1 & |i-j| = 1, i, j < 3,
		\\
		-2 & ij = 6,
		\\
		2 & i=j < 3,
		\\
		4 & i = j = 3,
		\\
		0 & \mathrm{otherwise.}
	\end{cases}
\end{equation}
The roots are ordered according to \eqref{ord} 
\eqn{ord21} \begin{matrix} &\b_{11}  &&&& \cr
	&\vee &&&&\cr
	&\b_{12} & > & \b_{22} &&\cr
	&\vee && \vee && \cr
	&\b_{13} & > & \b_{23} & > & \b_{33}=\g_3\cr
	&\vee &  & \vee  &  &    \cr
	&\a_{13} & > & \a_{23}=\g_2 &&\cr
	&\vee && &&\cr
	&\a_{12}=\g_1 &  &&  & \cr
\end{matrix}\ee
The long roots in terms of the simple roots are:
\eqn{longr} \b_{11} = 2\g_1+2\g_2+\g_3, ~ \b_{22}= 2\g_2+\g_3, ~\b_{33}= \g_3, \ee
the short:
\eqnn{shortr} &&\b_{12} = \g_1+2\g_2+\g_3, ~\b_{13} = \g_1+\g_2+\g_3, ~\b_{23} = \g_2+\g_3, \nn\\
&& \a_{12} = \g_1, ~ \a_{13} = \g_1+\g_2, ~\a_{23} = \g_2,
\eea

The basis of $ sp(3,\mathbb{C})$ in terms of those of $ sp(2,1)$ is given explicitly as follows:

\noindent
Cartan subalgebra
\begin{equation}
	H_1 = i M_2, \qquad H_2 = -i M_1 + \frac{1}{2} (A_0 - i M_2), \qquad
	H_3 = 2i M_1.
	\label{RelCartan}
\end{equation}

\noindent
Positive roots
\begin{align}
	E_{\gamma_1}  &= \frac{1}{2} (M_5 + i M_6), & 
	E_{\gamma_2}  &= -\frac{1}{2} (L^+_3-i L^+_4), & 
	E_{\gamma_3}  &= M_3 + i M_4,
	\nonumber \\
	E_{\beta_{11}} &= L^+_6 + i L^+_7, & 
	E_{\beta_{22}} &= L^+_6 - i L^+_7, & 
	E_{\beta_{12} } &= -i L^+_5,
	\nonumber \\
	E_{\beta_{13}} &= \frac{1}{2} (L^+_3 + i L^+_4), &
	E_{\beta_{23} } &= -\frac{1}{2}(L^+_1+iL^+_2), &
	E_{\alpha_{13} } &= - \frac{1}{2} (L^+_1- i L^+_2).
	\label{RelPositive}
\end{align}
Negative roots
\begin{align}
		E_{-\gamma_1 } &= -\frac{1}{2}(M_5-iM_6), & 
		E_{-\gamma_2 } &= \frac{1}{2} (L^-_3 + i L^-_4), & 
		E_{-\gamma_3 } &= -M_3 + i M_4,
		\nonumber \\
		E_{-\beta_{11} } &= -L^-_6 + i L^-_7, & 
		E_{-\beta_{22} } &= -L^-_6 - i L^-_7, & 
		E_{-\beta_{12} } &= -i L^-_5,
		\nonumber \\
		E_{-\beta_{13} } &= -\frac{1}{2} (L^-_3-iL^-_4), &
		E_{-\beta_{23} } &= \frac{1}{2} (L^-_1 - i L^-_2), &
		E_{-\alpha_{13} } &= \frac{1}{2} (L^-_1 + i L^-_2).
	\label{RelNegative}
\end{align}

From \eqref{RelCartan}, \eqref{RelPositive} and \eqref{RelNegative}, we see that
$
	{\cal M}_0^{\mathbb{C}} = sp(1,\mathbb{C}) \oplus sp(1,\mathbb{C})
$
is spanned by
\begin{equation}
	\{\  H_1, E_{\pm \gamma_1}  \}, \qquad \{\  H_3, E_{\pm \gamma_3}  \} \label{sp1component}
\end{equation}
and its triangular decomposition
$ {\cal M}_0^{\mathbb{C}} = ({\cal M}_{0}^+)^{\mathbb{C}} \oplus {\cal H}_0^{\mathbb{C}} \oplus ({\cal M}_{0}^-)^{\mathbb{C}} $ is readily understood.
Let ${\cal G}^{\mathbb{C}} = sp(3,\mathbb{C})$, then its triangular decomposition is given by
\begin{equation}
	{\cal G}_{\pm}^{\mathbb{C}} = ({\cal N}_0^{\pm})^{\mathbb{C}} \oplus ({\cal M}_{0}^{\pm})^{\mathbb{C}},
	\qquad
	{\cal H}^{\mathbb{C}} = {\cal H}_0^{\mathbb{C}} \oplus {\cal A}_0.
\end{equation}

 We return to the roots of ~$\cgc$~ which we need to divide them   into ~$\cm$-compact roots  and ~$\cm$-non-compact roots.
 The ~$\cm$-compact roots are the positive roots of ~$\cm^\bbc$, i.e.,
 \eqn{nonc}  \a_{12}=\g_1, ~\b_{33}= \g_3, \ee
 the rest are ~$\cm$-non-compact.

The HC parameters are:
\eqnn{hclabw}
&& m_{\b_{11}} = m_1+m_2+m_3 \equiv m_{13}, ~  \nn\\ &&
m_{\b_{22}} = m_2+m_3 \equiv m_{23}, ~
m_{\b_{33}} = m_3 \\
&& m_{\b_{12}} = m_1+2m_2+2m_3 \equiv m_{13,23}, ~  \nn\\ &&
m_{\b_{13}} = m_1+m_2+2m_3 \equiv m_{13,3}, ~  \nn\\
&& m_{\b_{23}} = m_2+2m_3 \equiv m_{23,3}, ~  \nn\\
&&m_{\a_{12}} = m_1 ,~
m_{\a_{23}} = m_2 , ~
m_{\a_{13}} = m_1+m_2 \equiv m_{12}
\eea

The main multiplets   are in 1-to-1 correspondence
with the finite-dimensional irreps of ~$sp(2,1)$, i.e., they are
labelled by  the  three  positive Dynkin labels     ~$m_1,m_2,m_3\,$.

We take ~$\chi_0 = \{ m_1,m_2,m_3 \}$. It has one embedded Verma module with HW ~$\L_a = \L_0-m_2\a_2$.
The number of ERs/GVMs   in a  main multiplet   is 12.
We give the whole multiplet as follows:

\eqnn{mainsp21z} && \chi_0 ~=~ \{ m_1,m_2,m_3 \} \\
&& \chi_a ~=~ \{ m_{12},-m_2,m_{23} \} \nn\\ &&
\chi_b ~=~ \{ m_{2},-m_{12},m_{13} \} \nn\\ &&
\chi_c ~=~ \{ m_{13,3},-m_{23,3},m_{23}\} \nn\\ &&
\chi_d ~=~ \{ m_{23,3},-m_{13,3},m_{13} \}  \nn\\ &&
\chi_e ~=~ \{ m_{13,23},-m_{23,3},m_{3} \} \nn\\ &&
\chi_f ~=~ \{ m_{23,3},-m_{13,23},m_{13} \} = \chh_d \nn\\ &&
\chi_g ~=~ \{ m_{13,23},-m_{13,3},m_{3} \} = \chh_e \nn\\ &&
\chi_h ~=~ \{ m_{2},-m_{13,23},m_{13} \} = \chh_b \nn\\ &&
\chi_i ~=~ \{ m_{13,3},-m_{13,23},m_{23} \} = \chh_c  \nn\\ &&
\chi_j ~=~ \{ m_{12},-m_{13,23},m_{23} \} = \chh_a \nn\\ &&
\chi_k ~=~ \{ m_{1},-m_{13,3},m_{3} \} = \chh_0 \nn
\eea
These multiplets  are presented in  Fig. 1.
On the figure each arrow represents an embedding between two Verma modules, ~$V^\L$ and ~$V^{\L'}$,~ the arrow pointing to the embedded module ~$V^{\L'}$.
Each arrow carries a number ~$n$, $n=1,2,3$, which indicates the level of the embedding,  ~$\L' = \L - m_n\,\b$. By our construction it also
represents the invariant differential operator ~$\cd_{n,\b}\,$, cf. \eqref{invop}.

Further, we note that there is an additional symmetry  w.r.t. to the dashed line in  Fig. 1.
It is relevant for the ERs and indicates the integral intertwining Knapp-Stein (KS) operators\cite{KnSt} acting between the
  spaces ~$\cc_{\chi^\mp}$~ in opposite directions:
\eqn{ksks} G^+_{KS} ~:~ \cc_{\chi^-} \llr \cc_{\chi^+}\ , \qquad
G^-_{KS} ~:~ \cc_{\chi^+} \llr \cc_{\chi^-} \ee
Note that the KS opposites are induced from the same irreps of ~$\cm$.

\begin{figure}[h]
	\centering
	\begin{tikzpicture}[scale=0.8]
		\coordinate (T) at (0,5);
		\coordinate (M) at (0,3);
		\coordinate (L) at (-3,1.5);  \coordinate (R) at (3,1.5);
		\coordinate (B) at (0,0);
		\coordinate (D) at (0,-5); 
		\coordinate (R1) at (6,0);
		\coordinate (R2) at (6,-2);
		\draw[dashed] (-4,-1) -- (7,-1);
		\draw[-Latex,thick] (T) node[above] {$\Lambda_0^-$} -- (M) node[right,yshift=5pt] {$\Lambda_a^-$} node[midway,right] {$2_2$};
		\draw[-Latex,thick] (M) -- (L) node[above left] {$\Lambda_b^-$} node[midway,above] {$1_{13}$};
		\draw[-Latex,thick] (M) -- (R) node[above right] {$\Lambda_c^-$} node[midway,above] {$3_{\widehat{22}}$};    
		\draw[-Latex,thick] (L) -- (B) node[above,yshift=5pt] {$\Lambda_d^-$} node[midway,above] {$3_{\widehat{22}}$};
		\draw[-Latex,thick] (R) -- (B) node[midway,above] {$1_{13}$};
		\draw[-Latex,thick] ($(T)+(D)$) -- ($(M)+(D)$) node[below,yshift=-3pt] {$\Lambda_d^+$} node[midway,right,yshift=7pt] {$2_{\widehat{12}}$};
		\draw[-Latex,thick] ($(M)+(D)$) -- ($(L)+(D)$) node[above left] {$\Lambda_b^+$} node[midway,above] {$3_{\widehat{11}}$};
		\draw[-Latex,thick] ($(M)+(D)$) -- ($(R)+(D)$) node[midway,above] {$1_{\widehat{23}}$} ;    
		\draw[-Latex,thick] ($(L)+(D)$) -- ($(B)+(D)$) node[above,yshift=5pt] {$\Lambda_a^+$} node[midway,above] {$1_{\widehat{23}}$};
		\draw[-Latex,thick] ($(R)+(D)$) -- ($(B)+(D)$) node[midway,above] {$3_{\widehat{11}}$};
		\draw[-Latex,thick] ($(T)+(0,-10)$) -- ($(M)+(0,-10)$) node[below] {$\Lambda_0^+$} node[midway,right] {$2_{\widehat{13}}$};
		\draw[-Latex,thick] (R) -- (R1) node[right] {$\Lambda_e^-$} node[midway,above] {$2_{\widehat{23}}$};
		\draw[-Latex,thick] (R1) -- (R2) node[right] {$\Lambda_e^+$} node[midway,right,yshift=7pt] {$1_{\widehat{12}}$};
		\draw[-Latex,thick] (R2) -- ($(R)+(D)$) node[above,yshift=5pt] {$\Lambda_c^+$} node[midway,above] {$2_{13}$};
	\end{tikzpicture}
	\caption{Main multiples for $sp(2,1).$}
\end{figure}
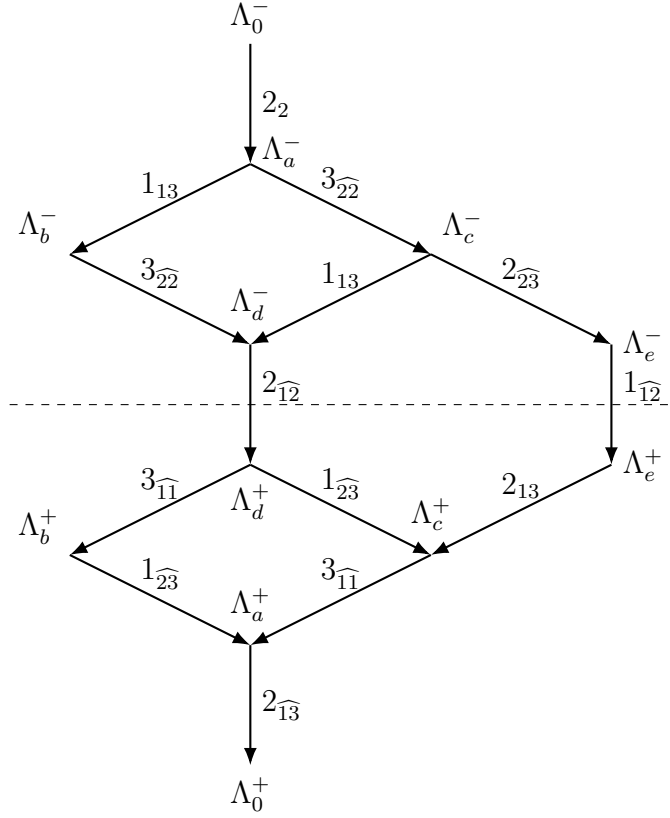

 We shall use also the other notation from \eqref{sgnd} for more compact exposition:
 \eqnn{mainsp21cz} && \chi^\pm_0 ~=~ [ m_1;\, \pm \ha m_{13,23}\,;m_3 ] \\
&& \chi^\pm_a ~=~ [ m_{12};\, \pm \ha m_{13,3}\,; m_{23} ] \nn\\ &&
\chi^\pm_b ~=~ [ m_{2};\, \pm \ha m_{23,3}\,; m_{13} ] \nn\\ &&
\chi^\pm_c ~=~ [ m_{13,3};\, \pm \ha m_{12}\,;  m_{23}] \nn\\ &&
\chi^\pm_d ~=~ [ m_{23,3};\, \pm \ha m_{2}\,; m_{13} ]  \nn\\ &&
\chi^\pm_e ~=~ [ m_{13,23};\, \pm \ha m_{1}\,; m_{3} ] \nn
\eea

{\it Remark:} ~Note that the pairs ~$\chi^\pm_d$~ and ~$\chi^\pm_e$~  are related by KS operators, but in each case the operator ~$G^+_{KS}$~   is degenerated into  a differential operator, namely, we have
  \eqna{multopjk} && \L^-_d ~~ {m_2\b_{12} \atop \longrightarrow } ~~\Lp_d \\ && \phantom{\L^-_j} \nn\\
&& \L^-_e ~~ {m_1\b_{12} \atop \longrightarrow } ~~\Lp_e \eena


\subsection{Reduced multiplets}

The reduced multiplet R21-1 is:

 \eqnn{mainsp21c1} && \chi^\pm_0 ~=~ [ 0;\, \pm \ha m_{23,23}\,;m_3 ] \\
&& \chi^\pm_a ~=~ [ m_{2};\, \pm \ha m_{23,3}\,; m_{23} ] = \chi^\pm_b \nn\\ &&
 \chi^\pm_d ~=~ [ m_{23,3};\, \pm \ha m_{2}\,; m_{23} ]  =\chi^\pm_c  \nn\\ &&
\chi_e ~=~ [ m_{23,23};\, 0\,; m_{3} ] \nn
\eea

\bigskip

The reduced multiplet R21-2 is:

\eqnn{mainsp21c2}
&& \chi^\pm_a ~=~ [ m_{1};\, \pm \ha m_{1,3,3}\,; m_{3} ] =  \chi^\pm_0 \\ &&
\chi^\pm_b ~=~ [ 0;\, \pm \ha m_{3,3}\,; m_{1,3} ] \nn\\ &&
\chi^\pm_c ~=~ [ m_{1,3,3};\, \pm \ha m_{1}\,;  m_{3}] =\chi^\pm_e \nn\\ &&
\chi_d ~=~ [ m_{3,3};\,  0\,; m_{1,3} ]  \nn
\eea
\bigskip

\begin{figure}[h]
	\centering
	\begin{tikzpicture}[scale=0.58]
		\foreach \y / \x in {15/2_2,12/3_{\widehat{22}},9/2_{\widehat{12}},6/3_{\widehat{11}},3/2_{\widehat{13}}}
		\draw[-Latex,thick] (0,\y) coordinate(\y) -- (0,\y-3) node[midway,right] {$\x$}; 
		\draw[-Latex,thick] (9) -- (3,7.5) coordinate(R) node[midway,above] {$2_{\widehat{23}}$};
		\draw[-Latex,thick] (R) -- (6) node[midway,below] {$2_{13}$};
		\node[left] at (15) {$\Lambda_0^-$};
		\node[left] at (12) {$\Lambda_a^-$};
		\node[left] at (9) {$\Lambda_d^-$};
		\node[left] at (6) {$\Lambda_d^+$};
		\node[left] at (3) {$\Lambda_a^+$};
		\node[left] at (0,0) {$\Lambda_0^+$};
		\node[right] at (R) {$\Lambda_e$}; 
		\begin{scope}[xshift=9.5cm,yshift=7cm]
			\coordinate (O) at (0,0);  \coordinate (A) at (3,1.5);
			\coordinate (B) at (3,-1.5); \coordinate (C) at (-3,1.5);
			\draw[-Latex,thick] (A) node[right] {$\Lambda_c^-$} -- (B) node[right] {$\Lambda_c^+$} node[midway,right] {$1_{\widehat{12}}$};
			\draw[-Latex,thick] (A) -- (O) node[above,yshift=+4pt] {$\Lambda_d$} node[midway,above] {$1_{13}$};
			\draw[-Latex,thick] ($(A)+(C)$) node[above] {$\Lambda_a^-$} -- ($(O)+(C)$) node[left] {$\Lambda_b^-$} node[midway,above] {$1_{13}$};
			\draw[-Latex,thick] (O) -- ($(A)+(-6,-3)$) node[left] {$\Lambda_b^+$} node[midway,above] {$3_{\widehat{11}}$};
			\draw[-Latex,thick] ($(A)+(0,-3)$) -- ($(O)+(0,-3)$) node[below] {$\Lambda_a^+$} node[midway,above] {$3_{\widehat{11}}$};
			\draw[-Latex,thick] (C) -- (O) node[midway,above] {$3_{\widehat{22}}$};
			\draw[-Latex,thick] ($(C)+(A)$) -- ($(O)+(A)$) node[midway,above] {$3_{\widehat{22}}$};
			\draw[-Latex,thick] ($(C)+(0,-3)$) -- ($(O)+(0,-3)$) node[midway,above] {$1_{\widehat{23}}$};
			\draw[-Latex,thick] (O) -- ($(A)+(0,-3)$) node[midway,above] {$1_{\widehat{23}}$};
		\end{scope}
		\begin{scope}[xshift=16.3cm,yshift=5.5cm]
			\coordinate (T) at (2.5,6) ;
			\coordinate (M1) at (2.5,4);
			\coordinate (OA) at (0,3); \coordinate (OD) at (5,3);
			\coordinate (O) at (0,0);  \coordinate (OR) at (5,0);
			\coordinate (M2) at (2.5,-1);
			\coordinate (B) at (2.5,-3);
			\draw[-Latex,thick] (T) node[above] {$\Lambda_0^-$} -- (M1) node[midway,right] {$2_2$};  
			\draw[-Latex,thick] (M1) node[below] {$\Lambda_c^-$} -- (OA) node[left] {$\Lambda_d^-$} node [midway, above] {$1_{13}$};
			\draw[-Latex,thick] (M1) -- (OD) node[right] {$\Lambda_e^-$} node[midway,above] {$2_{\widehat{23}}$};
			\draw[-Latex,thick] (OA) -- (O) node[left] {$\Lambda_d^+$} node[midway,right] {$2_{\widehat{12}}$};
			\draw[-Latex,thick] (OD) -- (OR) node[right] {$\Lambda_e^+$} node[midway,right] {$1_{\widehat{12}}$};
			\draw[-Latex,thick] (O) -- (M2) node[above] {$\Lambda_c^+$} node[midway,above] {$1_{\widehat{23}}$};
			\draw[-Latex,thick] (OR) -- (M2) node[midway,above] {$2_{13}$};
			\draw[-Latex,thick] (M2) -- (B) node[below] {$\Lambda_0^+$} node [midway,right] {$2_{\widehat{13}}$};
		\end{scope}
	\end{tikzpicture}
	\caption{Reduced multiplets R21-1 (left), R21-2 (middle) and R21-3 (right).}
	\label{Mpl2}
\end{figure}
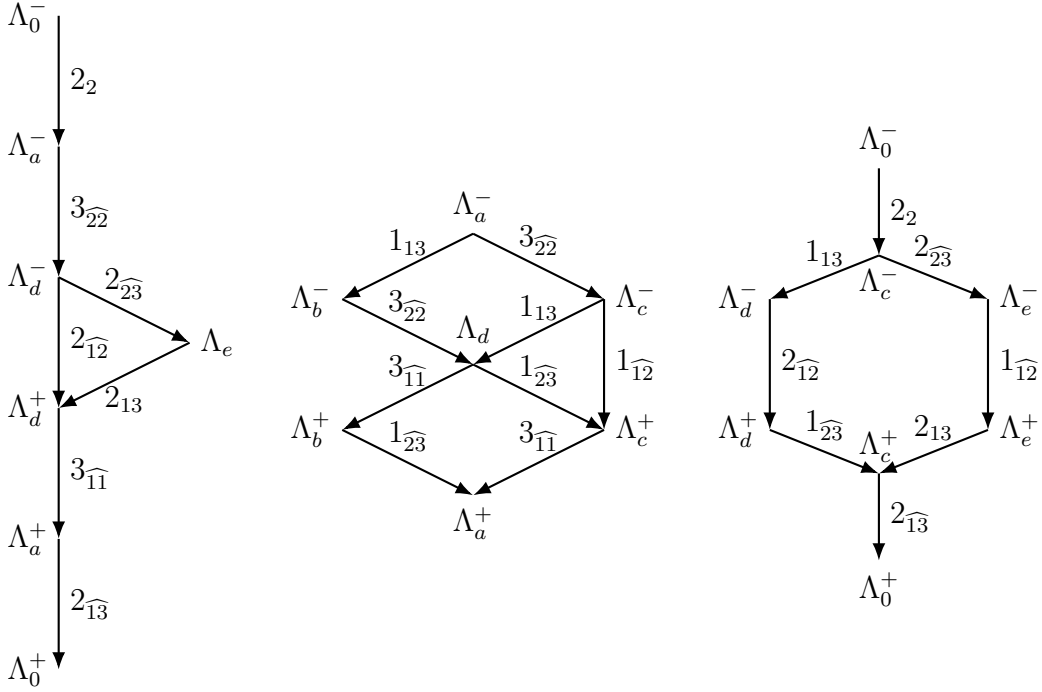

The reduced multiplet R21-3 is:

 \eqnn{mainsp21c3} && \chi^\pm_0 ~=~ [ m_1;\, \pm \ha m_{12,2}\,;0 ] \\  &&
\chi^\pm_c ~=~ [ m_{12};\, \pm \ha m_{12}\,;  m_{2}] =\chi^\pm_a  \nn\\ &&
\chi^\pm_d ~=~ [ m_{2};\, \pm \ha m_{2}\,; m_{12} ]  = \chi^\pm_b \nn\\ &&
\chi^\pm_e ~=~ [ m_{12,2};\, \pm \ha m_{1}\,; 0 ] \nn
\eea

\bigskip

The reduced multiplet R21-12 is:

  \eqnn{mainsp21c12}
&& \chi^\pm_a ~=~ [ 0 ;\, \pm \ha m_{3,3}\,; m_{3} ] = \chi^\pm_b =\chi^\pm_0 \\ &&
 \chi_d ~=~ [ m_{3,3};\,  0\,; m_{3} ]  =\chi^\pm_c  =\chi_e \nn
\eea

\bigskip

The reduced multiplet R21-13 is:
\eqnn{mainsp21c13} && \chi^\pm_0 ~=~ [ 0;\, \pm \ha m_{2,2}\,; 0  ] \\ &&
 \chi^\pm_d ~=~ [ m_{2};\, \pm \ha m_{2}\,; m_{2} ]  =\chi^\pm_c  =\chi^\pm_a \nn\\ &&
\chi_e ~=~ [ m_{2,2};\, 0\,; 0 ] \nn
\eea

 \bigskip

The reduced multiplet R21-23 is:
\eqnn{mainsp21c23}  &&
\chi^\pm_c ~=~ [ m_{1};\, \pm \ha m_{1}\,;  0] =\chi^\pm_e = \chi^\pm_a \\ &&
\chi_d ~=~ [ 0;\,  0\,; m_{1} ]  = \chi_b\nn
\eea

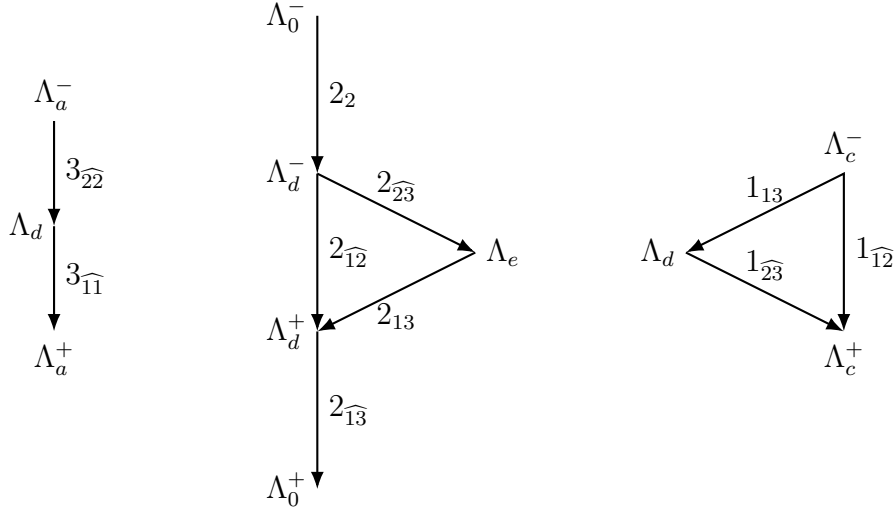
\begin{figure}[h]
	\centering
	\begin{tikzpicture}[scale=0.7]
	 \draw[-Latex,thick] (0,4) coordinate(T) -- (0,2) coordinate(M) node[midway,right]{$3_{\widehat{22}}$};
	 \draw[-Latex,thick] (M) -- (0,0) coordinate(B) node[midway,right] {$3_{\widehat{11}}$};
	 \node[above] at (T) {$\Lambda_a^-$};
	 \node[left] at (M) {$\Lambda_d$};
	 \node[below] at (B) {$\Lambda_a^+$};
		\begin{scope}[xshift=5cm,yshift=-3cm]
			\foreach \y / \x in {9/2_2,6/2_{\widehat{12}},3/2_{\widehat{13}}}
		\draw[-Latex,thick] (0,\y) coordinate(\y) -- (0,\y-3) node[midway,right] {$\x$};
		\draw[-Latex,thick] (6) -- (3,4.5) coordinate(R) node[midway,above] {$2_{\widehat{23}}$};
		\draw[-Latex,thick] (R) -- (3) node[midway,below] {$2_{13}$};
		\node[left] at (9) {$\Lambda_0^-$};
		\node[left] at (6) {$\Lambda_d^-$};
		\node[left] at (3) {$\Lambda_d^+$};
		\node[left] at (0,0) {$\Lambda_0^+$};
		\node[right] at (R) {$\Lambda_e$}; 
		\end{scope}
		\begin{scope}[xshift=15cm]
			\draw[Latex-Latex,thick] (-3,1.5) coordinate(L1) --(0,3) node[above](T){$\Lambda_c^-$} node[midway,above]{$1_{13}$} --(0,0) coordinate(L2) node[midway,right]{$1_{\widehat{12}}$};  
			\draw[-Latex,thick] (L1)--(L2) node[midway,above]{$1_{\widehat{23}}$};
			\node[below] at (L2) {$\Lambda_c^+$};
			\node[left] at (L1) {$\Lambda_d$};
		\end{scope}
	\end{tikzpicture}
	\caption{Reduced multiplets R21-12 (left), R21-13 (middle) and R21-23 (right).}
	\label{Mpl3}
\end{figure}

\bigskip
The corresponding diagrams are given in Fig.~\ref{Mpl2} and Fig.~\ref{Mpl3}.

\section{Explicit expressions of the \idos}

For the explicit expressions of the \idos\ we need first the explicit expressions of the singular vectors in terms of the root vectors.
For a class of roots these may be taken from \cite{Dobsin} adapting the notation.

\subsection{Right and left action}

We consider the reduced function $\varphi(g) $ defined on $g \in G^{\mathbb{C}}/B$ where $ G^{\mathbb{C}} = Sp(3,\mathbb{C}), B = \exp({\cal H}^{\mathbb{C}}) \exp({\cal G}_+^{\mathbb{C}}). $
We use the following parametrization
\begin{align}
	g &= \exp(Y_{\cal N}) \exp(Y_{\cal M}),
	\nonumber \\
	Y_{\cal N} &= x_1 E_{-\gamma_2} + x_2 E_{-\beta_{11}} + x_3 E_{-\beta_{22}} + x_4 E_{-\beta_{12}}
	\nonumber \\
	   & + x_5 E_{-\beta_{13}} + x_6 E_{-\beta_{23}} + x_7 E_{-\alpha_{13}},
	 \nonumber \\
	 Y_{\cal M} &= \xi_1 E_{-\gamma_1} + \xi_2 E_{-\gamma_3}.
	 \label{g-param}
\end{align}
The reason of this particular parametrization is
\begin{equation}
	Y_{\cal N}^k = 0 \quad (k \geq 3),  \quad Y_{\cal M}^k = 0 \quad (k \geq 2), \quad
	(Y_{\cal N} + Y_{\cal M})^k = 0 \quad (k \geq 6).
\end{equation}
This means that the parametrization \eqref{g-param} gives simpler matrix presentation than $ \exp(Y_{\cal N}+Y_{\cal M}). $
The explicit form of $ g$ is given by
\begin{eqnarray}
	g =
	\begin{pmatrix}
		1 & x_6 & 0 & -2\xi_2 & 0 & -x_5 -x_6 \xi_1
		\\
		0 & 1 & 0 & 0 & 0 & -\xi_1
		\\
		-x_7 & \chi & 1 & -x_5 + 2 x_7 \xi_2 & -\xi_1 & -2 x_2-\chi \xi_1
		\\
		0 & -x_1 & 0 & 1 & 0 & x_7 + x_1 \xi_1
		\\
		-x_1 & 2x_3 & 0 & -x_6 + 2x_1 \xi_2 & 1 & -2x_4-2 x_3 \xi_1+ \chi
		\\
		0 & 0 & 0 & 0 & 0 & 1
	\end{pmatrix}
\end{eqnarray}
where $ \chi = x_4 + \frac{1}{2} (x_1 x_5 - x_6 x_7). $

As mentioned in \S \ref{SEC:Pre}, due to the right covariance of $\varphi(g)$, we have
\begin{align}
   \widehat{H}_k \varphi &= h_k \varphi, \quad k = 1, 2, 3
   \nonumber \\
   \widehat{X} \varphi &= 0, \quad X \in {\cal G}_+^{\mathbb{C}}	
\end{align}
The right action of $ {\cal G}_-^{\mathbb{C}}$ on $ \varphi $ is calculated easily.
\begin{align}
	\widehat{E}_{-\gamma_1} &= \partial_{\xi_1},
	\nonumber \\
	\widehat{E}_{-\gamma_2} &= \partial_1 - \Big(\frac{1}{2} x_5 - x_7  \xi_2 \Big) \xi_1 \partial_2 + \Big( \frac{1}{2} x_6 - x_1 \xi_2 \Big) \partial_3
	\nonumber \\
	&+ \Big(\frac{1}{2} x_5 - x_7 \xi_2 - \big( \frac{1}{2}x_6-x_1 \xi_2 \big) \xi_1  \Big) \partial_4
	-2\xi_1 \xi_2 \partial_5 + 2\xi_2 \partial_6 - \xi_1 \partial_7,
	\nonumber \\
	\widehat{E}_{-\gamma_3} &= \partial_{\xi_2},
	\nonumber \\
	\widehat{E}_{-\beta_{11}} &= \partial_2,
	\nonumber \\
	\widehat{E}_{-\beta_{22}} &= \partial_3 + \xi_1^2 \partial_2 - 2\xi_1 \partial_4,
	\nonumber \\
	\widehat{E}_{-\beta_{12}} &= \partial_4 - \xi_1 \partial_2,
	\nonumber \\
	\widehat{E}_{-\beta_{13}} &= \partial_5 - \frac{1}{2} x_7 \partial_2 - \frac{1}{2}x_1 \partial_4,
	\nonumber \\
	\widehat{E}_{-\beta_{23}} &= \partial_6 + \frac{1}{2} x_7 \xi_1 \partial_2 - \frac{1}{2} x_1 \partial_3 -\frac{1}{2}(x_7-x_1 \xi_1) \partial_4 - \xi_1 \partial_5,
	\nonumber \\
	\widehat{E}_{-\alpha_{13}} &= \partial_7 + \Big( \frac{1}{2}x_5-x_7 \xi_2 \Big) \partial_2 + \Big( \frac{1}{2} x_6 - x_1 \xi_2 \Big) \partial_4 + 2 \xi_2 \partial_5
\end{align}
where $ \partial_k = \partial_{x_k}.$

We further make a change of variables:
\begin{align}
	y_1 &= x_1, & \qquad
	y_2 &= x_2 - \frac{1}{2} y_5  x_7,
	& \qquad
	y_3 &= x_3 - \frac{1}{2}x_1 y_6,
	\nonumber \\
	y_4 &= x_4, &
	y_5 &= x_5-2x_7 \xi_2,  &
	y_6 &= x_6 - 2 x_1 \xi_2,
	\nonumber \\
	y_7 &= x_7, &
	\eta_k &= \xi_k.  \label{NewVari}
\end{align}
Then, the right action has the simpler form
\begin{align}
	\widehat{E}_{-\gamma_1} &= \partial_{\eta_1},
	\nonumber \\
	\widehat{E}_{-\gamma_2} &= \partial_{y_1} + \frac{1}{2} (y_5-y_6\eta_1) \partial_{y_4} - \eta_1 \partial_{y_7},
	\nonumber \\
	\widehat{E}_{-\gamma_3} &= \partial_{\eta_2} + y_7^2 \partial_{y_2} + y_1^2 \partial_{y_3} - 2 y_7 \partial_{y_5} -2 y_1 \partial_{y_6},
	\nonumber \\
	\widehat{E}_{-\beta_{11}} &= \partial_{y_2},
	\nonumber \\
	\widehat{E}_{-\beta_{22}} &= \partial_{y_3} + \eta_1^2 \partial_{y_2} - 2\eta_1 \partial_{y_4},
	\nonumber \\
	\widehat{E}_{-\beta_{12}} &= \partial_{y_4} - \eta_1 \partial_{y_2},
	\nonumber \\
	\widehat{E}_{-\beta_{13}} &= \partial_{y_5} -y_7 \partial_{y_2} - \frac{y_1}{2} \partial_{y_4},
	\nonumber \\
	\widehat{E}_{-\beta_{23}} &= \partial_{y_6} + y_7 \eta_1 \partial_{y_2} - y_1 \partial_{y_3} - \frac{1}{2} (y_7-y_1 \eta_1) \partial_{y_4},
	\nonumber \\
	\widehat{E}_{-\alpha_{13}} &= \partial_{y_7} + \frac{y_6}{2} \partial_{y_4}.
	\label{RighActionSimple}
\end{align}

Next we need the left action of $\mathcal{G}^{\mathbb{C}}$, denoted by $\pi_L,$ which is the infinitesimal version of \eqref{lrr}.

\bigskip\noindent
$-$ Cartan subalgebra
\begin{align}
	\pi_L(H_1) &= -x_1 \partial_1 + 2x_2 \partial_2 - 2 x_3 \partial_3 + x_5 \partial_5-x_6 \partial_6 + x_7 \partial_7 + 2\xi_1 \partial_{\xi_1}-h_1,
	\nonumber \\
	\pi_L(H_2) &=2x_1 \partial_1 + 2x_3 \partial_3 + x_4 \partial_4 - x_5 \partial_5 + x_7 \partial_7 -\xi_1 \partial_{\xi_1} - 2 \xi_2 \partial_{\xi_2} - h_2,
	\nonumber \\
	\pi_L(H_3) &= 2(-x_1\partial_1 + x_5\partial_5+ x_6\partial_6 -x_7\partial_7 + 2 \xi_2 \partial_{\xi_2}) - h_3.
	\label{LeftCartan}
\end{align}
In the new variables \eqref{NewVari}, these left action preserve their form.
Namely, one may replace $ (x_k, \xi_k) \to (y_k, \eta_k) $ in \eqref{LeftCartan}.

\medskip\noindent
$-$ Positive roots
\begin{align}
	\pi_L(E_{\gamma_1}) &= x_7 \partial_1 + x_4 \partial_3 + 2x_2 \partial_4 + x_5 \partial_6 + \xi_1^2 \partial_{\xi_1} - \xi_1 h_1,
	\nonumber \\
	\pi_L(E_{\gamma_2}) &= x_1 \pi_L(H_2) -x_1 \Big( x_1 \partial_{x_1} + x_3 \partial_{x_3} + \frac{1}{2} (x_4 \partial_{x_4} - x_5 \partial_{x_5}) -x_6 \partial_{x_6} \Big)
	\nonumber \\
	&+ \frac{x_7}{4} (2x_4 + x_1 x_5 - x_6 x_7) \partial_2
	 + \Big( x_3 x_7 + \frac{x_1}{4}(x_1 x_5 - x_6 x_7) \Big) \partial_4
	 \nonumber \\
	 &+ \Big( x_4 + \frac{3}{2} x_6 x_7 \Big) \partial_5 + 2x_3 \partial_6
	  - x_7 \partial_{\xi_1} + x_6 \partial_{\xi_2},
	\nonumber \\
	\pi_L(E_{\gamma_3}) &= -2 x_6 \partial_1 - 2x_5 \partial_7 + 4\xi_2^2 \partial_{\xi_2} - 2\xi_2 h_3.
\end{align}
The change of  variables to \eqref{NewVari} gives the following expressions of the left action:
\begin{align}
	\pi_L(E_{\gamma_1}) &= y_7 \partial_{y_1} + \Big( y_4 -\frac{1}{2}(y_1 y_5 + y_6 y_7)  \Big) \partial_{y_3} + (2y_2 + y_5 y_7) \partial_{y_4}
	\nonumber \\
	& + y_5 \partial_{y_6} + \eta_1^2 \partial_{\eta_1} - \eta_1 h_1,
	\nonumber \\
	\pi_L(E_{\gamma_2}) &= y_1 \pi_L(H_2) - y_1 \Big(y_1 \partial_{y_1} + 2y_3 \partial_{y_3}  +\frac{1}{2}(y_4 \partial_{y_4}-y_5 \partial_{y_5}) -2\eta_2 \partial_{\eta_2} )
	\nonumber \\
	&+ \Big( y_3 y_7 + \frac{y_1}{4}(y_1 y_5 + y_6 y_7) \Big) \partial_{y_4}
	 + \Big( y_4 - \frac{1}{2} y_6 y_7 \Big) \partial_{y_5}
	 \nonumber \\
	& + 2y_3 \partial_{y_6} - y_7 \partial_{\eta_1} + y_6 \partial_{\eta_2},
	\nonumber \\
	\pi_L(E_{\gamma_3}) &= 2\eta_2 \pi_L(H_3) - 2 y_6 \partial_{y_1} + y_5^2 \partial_{y_2} +y_6^2 \partial_{y_3} -2 y_5 \partial_{y_7} - 4\eta_2^2 \partial_{\eta_2}.
\end{align}

\medskip\noindent
$-$ Negative roots
\begin{align}
	\pi_L(E_{-\gamma_1}) &= -\partial_{\xi_1} + x_4 \partial_{2} + 2 x_3 \partial_{4} + x_6 \partial_{5} + x_1 \partial_{7},
	\nonumber \\
	\pi_L(E_{-\gamma_2}) &= -\partial_{1} + \frac{1}{2} x_5 \partial_{4} + \frac{1}{2} x_6 \partial_{3},
	\nonumber \\
	\pi_L(E_{-\gamma_3}) &= -\partial_{\xi_2} - 2 x_7 \partial_{5} - 2x_1 \partial_{6},
	\nonumber \\
	\pi_L(E_{-\beta_{11}}) &= -[\pi_L(E_{-\gamma_1}), \pi_L(E_{-\beta_{12}}) ]  =-\partial_2,
	\nonumber \\
	\pi_L(E_{-\beta_{22}}) &= -[\pi_L(E_{-\gamma_2}), \pi_L(E_{-\beta_{23}}) ] = -\partial_3,
	\nonumber \\
	\pi_L(E_{-\beta_{12}}) &= -[\pi_L(E_{-\gamma_1}), \pi_L(E_{-\beta_{13}}) ]  =-\partial_4,
	\nonumber \\
	\pi_L(E_{-\beta_{13}}) &= \frac{1}{2} [\pi_L(E_{-\gamma_3}), \pi_L(E_{-\alpha_{13}})]
	    = -\partial_5 - \frac{1}{2} x_7 \partial_2 - \frac{1}{2}x_1 \partial_4,
	\nonumber \\
    \pi_L(E_{-\beta_{23}}) &= -\frac{1}{2} [\pi_L(E_{-\gamma_2}), \pi_L(E_{-\gamma_3})]
        = -\partial_6 - \frac{1}{2} x_1 \partial_3 - \frac{1}{2}x_7 \partial_4,
	\nonumber \\
    \pi_L(E_{-\alpha_{13}}) &= -[\pi_L(E_{-\gamma_1}), \pi_L(E_{-\gamma_2})]
    = -\partial_7 + \frac{1}{2} x_6 \partial_4 + \frac{1}{2}x_5 \partial_2.
\end{align}
In terms of the variables \eqref{NewVari}
\begin{align}
	\pi_L(E_{-\gamma_1}) &= -\partial_{\eta_1} + \Big( y_4-\frac{1}{2}(y_1 y_5 +y_6 y_7)  \Big) \partial_{y_2} + (2y_3+y_1 y_6) \partial_{y_4} + y_6 \partial_{y_5} + y_1 \partial_{y_7},
	\nonumber \\
	\pi_L(E_{-\gamma_2}) &= -\partial_{y_1} + y_6 \partial_{y_3} + \frac{1}{2}(y_5 + 2y_7 \eta_2) \partial_{y_4} + 2\eta_2 \partial_{y_6},
	\nonumber \\
	\pi_L(E_{-\gamma_3}) &= -\partial_{\eta_2},
	\nonumber \\
	\pi_L(E_{-\beta_{11}}) &=-\partial_{y_2},
	\nonumber \\
	\pi_L(E_{-\beta_{22}}) &= -\partial_{y_3},
	\nonumber \\
	\pi_L(E_{-\beta_{12}}) &=-\partial_{y_4},
	\nonumber \\
	\pi_L(E_{-\beta_{13}}) &= -\partial_{y_5} - \frac{y_1}{2} \partial_{y_4},
	\nonumber \\
	\pi_L(E_{-\beta_{23}}) &= -\partial_{y_6} - \frac{y_7}{2} \partial_{y_4},
	\nonumber \\
	\pi_L(E_{-\alpha_{13}}) &= -\partial_{y_7} + y_5 \partial_{y_2} + \frac{1}{2}(y_6+2y_1 \eta_2) \partial_{y_4} +2\eta_2 \partial_{y_5}.
\end{align}

\medskip
Before closing this subsection, we give the finite dimensional representations of $sp(1)$ in \eqref{sp1component}. 
The HWs and the Dynkin labels are related as follows:
\begin{equation}
	(h_1, h_2, h_3)  = (m_1-1, m_2-1,  2(m_3-1))
\end{equation}	
with $ h_j = \Lambda(H_j). $
The equivalent expression is
\begin{align}
	\Lambda &= \sum_{j=1}^3 c_j \gamma_j,
	\nonumber \\
	c_1 &= m_1 + m_2 + m_3-3,
	\nonumber \\
	c_2 &= m_1 + 2 m_2 + 2m_3 -5,
	\nonumber\\
	c_3 &= m_2 +\frac{1}{2}(m_1+3m_3)-3.
\end{align}

Repeated action of $ \pi_L(E_{\gamma_1}) $ on a constant function $C$ gives a finite dimensional representation of $ sp(1)$:
For $ h_1 \geq 0 \;(m_1 \geq 1)$
\begin{equation}
	\big( \pi_L(E_{\gamma_1}) \big)^k C =
	\begin{cases}
		\displaystyle (-\xi_1)^k \frac{h_1!}{(h_1-k)!}C,
		& k \leq h_1
		\\[15pt]
		0 & k > h_1
	\end{cases}
\end{equation}
We now rescale the elements
\begin{equation}
	\tilde{H}_3 = \frac{1}{2}, \qquad \tilde{E}_{\pm \gamma_3} = \frac{1}{2} E_{\pm \gamma_3}.
\end{equation}
Then, they satisfy
\begin{equation}
	[\tilde{H}_3, \tilde{E}_{\pm \gamma_3} ] = \pm 2 E_{\pm \gamma_3},
	\qquad
	[\tilde{E}_{\gamma_3}, \tilde{E}_{-\gamma_3} ] = \tilde{H}_3
\end{equation}
and
\begin{equation}
	 \tilde{H}_3 \varphi = \tilde{h}_3 \varphi, \qquad \tilde{h}_3 = \frac{1}{2} h_3=m_3-1.
\end{equation}
It follows that we have a finite dimensional representation of $sp(1)$ for $ \tilde{h}_3 \geq 0 \; (m_3 \geq 1)$
\begin{equation}
		\big( \pi_L(E_{\gamma_3}) \big)^k C =
		\begin{cases}
			\displaystyle (-2\xi_2)^k \frac{\tilde{h}_3!}{(\tilde{h}_3-k)!}C & k \leq \tilde{h}_3,
			\\[15pt]
			0 & k > \tilde{h}_3
		\end{cases}		
\end{equation}

We remark that
\begin{equation}
	H_i, \quad \tilde{H}_3, \quad E_{\pm \gamma_1}, \quad  \tilde{E}_{\pm \gamma_3}, \quad i = 1, 2
	\label{Chev}
\end{equation}
provides the Chevalley basis of $ sp(3,\mathbb{C}). $

%
\subsection{Invariant differential operators (IDOs)}

We present singular vectors  in terms of the Chevalley basis \eqref{Chev} as in \cite{Dobsin,VKD1}.
Then present the corresponding IDO for each positive root. 
In the following, $s_k(\beta)$ denotes the Weyl reflection with respect to the simple root $\gamma_k$.

\medskip\noindent
i) for the simple roots, the SVs are specified by the Dynkin label
\begin{equation}
	v^s_{m_i, \gamma_i}= (E_{-\gamma_i})^{m_i}v_0, \qquad v^s_{m_3, \gamma_3} = (\tilde{E}_{-\gamma_3})^{m_3} v_0, \quad i = 1,2
\end{equation}
The corresponding IDOs are given as follows:
\begin{align}
	{\cal D}_{m_1,\gamma_1} &= \partial_{\eta_1}^{m_1},
	\nonumber\\
    {\cal D}_{m_2,\gamma_2} &=\big(\partial_{y_1} + \frac{1}{2} (y_5-y_6\eta_1) \partial_{y_4} - \eta_1 \partial_{y_7} \big)^{m_2},
    \nonumber\\
    {\cal D}_{m_3,\gamma_3} &= (\partial_{\eta_2} + y_7^2 \partial_{y_2} + y_1^2 \partial_{y_3} - 2 y_7 \partial_{y_5} -2 y_1 \partial_{y_6})^{m_3}.
\end{align}

\noindent
ii) $ \alpha_{13} = \gamma_1 + \gamma_2 = s_2(\gamma_1) $

Writing the HC parameter $ m \equiv m_{12} = m_1 + m_2, $ then SV is given by
\begin{align}
	v^s_{m,\alpha_{13}} &= \sum_{k=0}^m a_k (E_{-\gamma_2})^{m-k} (E_{-\gamma_1})^m (E_{-\gamma_2})^{k} v_0,
	\nonumber\\
	a_k &= (-1)^k
	\begin{pmatrix}
		m \\ k
	\end{pmatrix}
    \frac{(\Lambda+\rho)(H_2)}{(\Lambda+\rho)(H_2)-k}
    = (-1)^k
	\begin{pmatrix}
	m \\ k
\end{pmatrix}
\frac{m_2}{m_2-k}
\end{align}
with $ \rho = 3\gamma_1 + 5 \gamma_2 + 3 \gamma_3. $ 
While the simpler form of IDO is obtained by the alternative presentation $ m = s_1(\gamma_2) $ for which the equivalent SV read as follows:
\begin{align}
	v^s_{m,\alpha_{13}} &= \sum_{k=0}^m a_k (E_{-\gamma_1})^{m-k} (E_{-\gamma_2})^m (E_{-\gamma_1})^{k} v_0,
	\nonumber\\
	a_k &= (-1)^k
	\begin{pmatrix}
		m \\ k
	\end{pmatrix}
	\frac{(\Lambda+\rho)(H_1)}{(\Lambda+\rho)(H_1)-k}
	= (-1)^k
	\begin{pmatrix}
		m \\ k
	\end{pmatrix}
	\frac{m_1}{m_1-k}.
\end{align}

The IDO obtained from this SV is given by
\begin{align}
	{\cal D}_{m,\alpha_{13}} &=  \sum_{k=0}^m a_k \partial_{\eta_1}^{m-k} (\partial_{y_1} + \frac{1}{2} (y_5-y_6\eta_1) \partial_{y_4} - \eta_1 \partial_{y_7})^m \partial_{\eta_1}^{k}.
\end{align}

\noindent
iii) $ \beta_{23} = \gamma_2 + \gamma_3 = s_3(\gamma_2), $ $ m \equiv m_{23,3} = m_2 + 2 m_3$
\begin{align}
	v^s_{m,\beta_{23}} &= \sum_{k=0}^m a_k (\tilde{E}_{-\gamma_3})^{m-k} (E_{-\gamma_2})^{m} (\tilde{E}_{-\gamma_3})^{k} v_0,
	\nonumber\\
	a_k &= (-1)^k
	\begin{pmatrix}
		m \\ k
	\end{pmatrix}
	\frac{(\Lambda+\rho)(\tilde{H}_3)}{(\Lambda+\rho)(\tilde{H}_3)-k}
	= (-1)^k
	\begin{pmatrix}
		m \\ k
	\end{pmatrix}
	\frac{m_3}{m_3-k}.
\end{align}
\begin{align}
	{\cal D}_{m,\beta_{23}} &= \sum_{k=0}^m a_k (\partial_{\eta_2} + y_7^2 \partial_{y_2} + y_1^2 \partial_{y_3} - 2 y_7 \partial_{y_5} -2 y_1 \partial_{y_6})^{m-k}
	\nonumber \\
	& \times \big(\partial_{y_1} + \frac{1}{2} (y_5-y_6\eta_1) \partial_{y_4} - \eta_1 \partial_{y_7} \big)^{m}
	\nonumber \\
	& \times  (\partial_{\eta_2} + y_7^2 \partial_{y_2} + y_1^2 \partial_{y_3} - 2 y_7 \partial_{y_5} -2 y_1 \partial_{y_6})^{k}.
\end{align}

\noindent
iv) $ \beta_{22} = 2\gamma_2 + \gamma_3 = s_2(\gamma_3), $ $ m \equiv m_{23} = m_2 + m_3$
\begin{align}
	v^s_{m,\beta_{22}} &= \sum_{k=0}^{2m} c_k (E_{-\gamma_2})^{2m-k} (\tilde{E}_{-\gamma_3})^{m} (E_{-\gamma_2})^{k} v_0,
	\nonumber\\
	c_k &= (-1)^k
	\begin{pmatrix}
		2m \\ k
	\end{pmatrix}
	\frac{(\Lambda+\rho)(H_2)}{(\Lambda+\rho)(H_2)-k}
	= (-1)^k
	\begin{pmatrix}
		2m \\ k
	\end{pmatrix}
	\frac{m_2}{m_2-k}.
\end{align}
\begin{align}
	{\cal D}_{m,\beta_{22}} &= \sum_{k=0}^{2m} c_k
	\big(\partial_{y_1} + \frac{1}{2} (y_5-y_6\eta_1) \partial_{y_4} - \eta_1 \partial_{y_7} \big)^{2m-k}
	\nonumber \\
	& \times (\partial_{\eta_2} + y_7^2 \partial_{y_2} + y_1^2 \partial_{y_3} - 2 y_7 \partial_{y_5} -2 y_1 \partial_{y_6})^{m}
	\nonumber \\
	& \times
	\big(\partial_{y_1} + \frac{1}{2} (y_5-y_6\eta_1) \partial_{y_4} - \eta_1 \partial_{y_7} \big)^{k}.
\end{align}

\noindent
v) $ \beta_{13} = \gamma_1 + \gamma_2 + \gamma_3 = s_3 s_2 (\gamma_1),$
$ m \equiv m_{13,3} = m_1 + m_2 + 2m_3 $
\begin{align}
	v^s_{m,\beta_{13}} &= \sum_{k_1, k_2=0}^m a_{k_1,k_2}
	(\tilde{E}_{-\gamma_3})^{m-k_1} (E_{-\gamma_2})^{m-k_2} (E_{-\gamma_1})^{m}(E_{-\gamma_2})^{k_2} (\tilde{E}_{-\gamma_3})^{k_1} v_0,
	\nonumber \\
	a_{k_1,k_2} &= (-1)^{k_1+k_2}
	\begin{pmatrix}
		m \\ k_1
	\end{pmatrix}
	\begin{pmatrix}
	    m \\ k_2
    \end{pmatrix}
     \frac{(\Lambda+\rho)(\tilde{H}_3)}{(\Lambda+\rho)(\tilde{H}_3)-k_1}
     \frac{(\Lambda+\rho)(H_2+2\tilde{H}_3)}{(\Lambda+\rho)(H_2+2\tilde{H}_3)-k_2}
     \nonumber \\
     &= (-1)^{k_1+k_2}
	\begin{pmatrix}
		m \\ k_1
	\end{pmatrix}
	\begin{pmatrix}
		m \\ k_2
	\end{pmatrix}
	\frac{m_3}{m_3-k_1}
	\frac{m_2 + 2m_3}{m_2 + 2m_3-k_2}.
\end{align}
\begin{align}
	{\cal D}_{m,\beta_{13}} &=  \sum_{k_1, k_2=0}^m a_{k_1,k_2}
	(\partial_{\eta_2} + y_7^2 \partial_{y_2} + y_1^2 \partial_{y_3} - 2 y_7 \partial_{y_5} -2 y_1 \partial_{y_6})^{m-k_1}
	\nonumber \\
	& \times \big(\partial_{y_1} + \frac{1}{2} (y_5-y_6\eta_1) \partial_{y_4} - \eta_1 \partial_{y_7} \big)^{m-k_2} \partial_{\eta_1}^{m}
	\nonumber \\
	& \times
	\big(\partial_{y_1} + \frac{1}{2} (y_5-y_6\eta_1) \partial_{y_4} - \eta_1 \partial_{y_7} \big)^{k_2}
	\nonumber \\
	& \times
	(\partial_{\eta_2} + y_7^2 \partial_{y_2} + y_1^2 \partial_{y_3} - 2 y_7 \partial_{y_5} -2 y_1 \partial_{y_6})^{k_1}.
\end{align}

\noindent
vi) $ \beta_{11} = 2\gamma_1 + 2 \gamma_2 + \gamma_3, $  $ m \equiv m_{m_{13}} = m_1 + m_2 + m_3 $
\begin{align}
	v^s_{m,\beta_{11}} &= \sum_{k_1, k_2=0}^{2m} c_{k_1,k_2}
	(E_{-\gamma_1})^{2m-k_1} (E_{-\gamma_2})^{2m-k_2} (\tilde{E}_{-\gamma_3})^{m} (E_{-\gamma_2})^{k_2} (E_{-\gamma_1})^{k_1} v_0,
	\nonumber \\
	c_{k_1,k_2} &= (-1)^{k_1+k_2}
	\begin{pmatrix}
		2m \\ k_1
	\end{pmatrix}
	\begin{pmatrix}
		2m \\ k_2
	\end{pmatrix}
	\frac{(\Lambda+\rho)(H_1)}{(\Lambda+\rho)(H_1)-k_1}
	\frac{(\Lambda+\rho)(H_1+H_2)}{(\Lambda+\rho)(H_1+H_2)-k_2}
	\nonumber \\
	&= (-1)^{k_1+k_2}
	\begin{pmatrix}
		2m \\ k_1
	\end{pmatrix}
	\begin{pmatrix}
		2m \\ k_2
	\end{pmatrix}
	\frac{m_1}{m_1-k_1}
	\frac{m_1 + m_2}{m_1+m_2-k_2}.
\end{align}
\begin{align}
	{\cal D}_{m,\beta_{11}} &= \sum_{k_1, k_2=0}^{2m} c_{k_1,k_2} \,
		\partial_{\eta_1}^{2m-k_1}
		\big(\partial_{y_1} + \frac{1}{2} (y_5-y_6\eta_1) \partial_{y_4} - \eta_1 \partial_{y_7} \big)^{2m-k_2}
		\nonumber \\
		& \times
		(\partial_{\eta_2} + y_7^2 \partial_{y_2} + y_1^2 \partial_{y_3} - 2 y_7 \partial_{y_5} -2 y_1)^m
        \nonumber \\
        & \times
       \big(\partial_{y_1} + \frac{1}{2} (y_5-y_6\eta_1) \partial_{y_4} - \eta_1 \partial_{y_7} \big)^{k_2}
        \partial_{\eta_1}^{k_1}.
\end{align}

\begin{align}
	\widehat{E}_{-\gamma_1} &= \partial_{\eta_1},
	\nonumber \\
	\widehat{E}_{-\gamma_2} &= \partial_{y_1} + \frac{1}{2} (y_5-y_6\eta_1) \partial_{y_4} - \eta_1 \partial_{y_7},
	\nonumber \\
	\widehat{E}_{-\gamma_3} &= \partial_{\eta_2} + y_7^2 \partial_{y_2} + y_1^2 \partial_{y_3} - 2 y_7 \partial_{y_5} -2 y_1 \partial_{y_6},
\end{align}

\section{Conclusion}

We presented explicit expressions for IDOs of $Sp(2,1)$ which is the first nontrivial member of split rank one non-compact group $ Sp(n,1). $ 
This is done by obtaining explicit expressions of singular vectors in the reducible Verma modules over the complexification of the Lie algebra $sp(3,\mathbb{C}). $ 
One may repeat this systematic construction of IDOs for $ n > 2. $ 
Obviously, it will become harder to find the explicit expression of singular vectors as $n$ goes larger. 
The present work will provide some guidance for the case of larger values of $n.$  

	Finally, we comment that IDOs are also arising from subsingular vectors and this is an interesting and challenging future work. The Weyl group $W_3$ of $sp(3,\mathbb{C})$ has 48 elements and one has to  find all pairs of elements of $W_3$  with nontrivial  Kazhdan-Lusztig polynomials taking into account that not all such pairs give rise to subsingular vectors (see \cite{Subsing} for detail).


\bigskip\bigskip

\noindent {\bf Acknowledgments.}

N. A. is supported by JSPS KAKENHI Grant Number JP23K03217.


\bigskip\bigskip

\begin{thebibliography}{99.}


\bibitem{DMPPT} V.K. Dobrev, G. Mack, V.B. Petkova, S.G. Petrova and I.T.
Todorov,
{\it Harmonic Analysis on the n-Dimensional Lorentz Group
and Its Applications to Conformal Quantum Field Theory}, Lecture
Notes in Physics, No 63, 280 pages (Springer
Verlag, Berlin-Heidelberg-New York, 1977);\\
 V.K. Dobrev and V.B. Petkova, 
Rep. Math. Phys. {\bf 13} (1978) 233-277.

 \bibitem{Zhel} D.P. Zhelobenko, Discrete symmetry operators for reductive Lie groups, Izv.
Akad. Nauk SSSR Ser. Mat. {\bf 40} (1976) 1055-1083 (in Russian).

\bibitem{Dobf4}V.K. Dobrev, Invariant Differential Operators for the Real Exceptional Lie Algebra $F''_4$,
arXiv:2109.08395, Contribution to  Peter Suranyi 87th Birthday Festschrift:
"A Life in Quantum Field Theory", Eds. R. Wijewardhana, Ph. Argyres, G. Dunne, G.  Semenoff. 
 (World Scientific, 2023),
Pages: 149–166, (Chapter 9). 

\bibitem{Pajas}P. Pajas, 
J. Math. Phys. {\bf 10} (1969) 1777–1788. 

\bibitem{SUGRA}H.-Y. Chang, E. Sezgin and Y. Tanii, 
J. High Energ. Phys. {\bf 2023}  (2023) 172.

\bibitem{Dob}V.K. Dobrev,
Rept. Math. Phys. {\bf 25}  (1988) 159-181; first as ICTP Trieste
preprint IC/86/393 (1986).

\bibitem{Lan}R.P. Langlands, {\it On the classification of irreducible
representations of real algebraic groups}, Math. Surveys and
Monographs, Vol. 31 (AMS, 1988), first as IAS Princeton preprint
(1973).

\bibitem{Zhea}D.P. Zhelobenko, {\it Harmonic Analysis on Semisimple
Complex Lie Groups}, (Moscow, Nauka, 1974, in Russian).

\bibitem{KnZu}A.W. Knapp and G.J. Zuckerman, ``Classification theorems
for representations of semisimple groups'',
 in: Lecture Notes in Math., Vol. 587 (Springer, Berlin,
1977) pp. 138-159; ~``Classification of irreducible tempered
representations of semisimple groups'', Ann. Math. {\bf 116} (1982)
389-501.

\bibitem{Har}Harish-Chandra, ``Discrete series for semisimple Lie
groups: II'', Ann. Math. {\bf 116} (1966) 1-111.

\bibitem{EHW}T. Enright, R. Howe and W. Wallach, "A classification of
unitary highest weight modules", in: {\it Representations of
Reductive Groups}, ed. P. Trombi (Birkh\"auser, Boston, 1983) pp.
97-143.



\bibitem{Knapp}A.W. Knapp, {\it Representation Theory of Semisimple
Groups (An Overview Based on Examples)}, (Princeton Univ. Press,
1986).


\bibitem{Dobmul}V.K. Dobrev, Lett. Math. Phys. {\bf 9}
(1985) 205-211; J. Math. Phys. {\bf 26} (1985) 235-251.


\bibitem{BGG}I.N. Bernstein, I.M. Gel'fand and S.I. Gel'fand,
``Structure of representations generated by highest weight
vectors'', Funkts. Anal. Prilozh. {\bf 5} (1) (1971) 1-9; English
translation: Funct. Anal. Appl. {\bf 5} (1971) 1-8.

\bibitem{Dix}J. Dixmier, {\it Enveloping Algebras}, (North Holland, New
York, 1977).

\bibitem{Sata}I. Satake,
Ann. Math. {\bf 71} (1960) 77-110.


\bibitem{VKD1} Vladimir K. Dobrev, {\it Invariant Differential Operators,
Volume 1: Noncompact Semisimple Lie Algebras and Groups}, De Gruyter
Studies in Mathematical Physics vol. 35 (De Gruyter, Berlin, Boston,
2016,
ISBN 978-3-11-042764-6), 408 + xii pages.%


\bibitem{KnSt}A.W. Knapp and E.M. Stein, ``Intertwining operators for
semisimple groups'', Ann. Math. {\bf 93} (1971) 489-578; II : Inv.
Math. {\bf 60} (1980) 9-84.

\bibitem{Dobsin}V.K. Dobrev,
Lett. Math. Phys. {\bf 22}   (1991) 251-266.


\bibitem{Subsing}V.K. Dobrev, 
 ``Kazhdan-Lusztig Polynomials, Subsingular Vectors and Conditionally Invariant (q-Deformed) Equations", In: Gruber, B., Ramek, M. (eds) Symmetries in Science IX. Springer, Boston, MA, 1997.
\end{thebibliography}
\end{document}